\pdfoutput=1
\RequirePackage{ifpdf}
\ifpdf % We are running pdfTeX in pdf mode
\documentclass[pdftex]{sigma}
\else
\documentclass{sigma}
\fi

\numberwithin{equation}{section}

\newtheorem{Theorem}{Theorem}[section]
\newtheorem{conj}[Theorem]{Conjecture}
\newtheorem{Questions}[Theorem]{Questions}

\newcommand{\hyp}[5]{{}_{#1}F_{#2} \left(\begin{array}{@{}c@{}} {#3}
\\[0.2cm] {#4} \end{array};{#5}\right)}
\newcommand{\choo}[2]{\raisebox{3pt}{$\left\{\substack{#1\\#2}\right\}$}}
\newcommand{\qhyp}[6]{{}_{#1}\phi_{#2} \left(\begin{array}{@{}c@{}} {#3}
\\[0.2cm] {#4} \end{array};{#5},{#6}\right)}

\begin{document}

\allowdisplaybreaks

\renewcommand{\thefootnote}{$\star$}

\renewcommand{\PaperNumber}{071}

\FirstPageHeading

\ShortArticleName{Report from the Open Problems Session at OPSFA13}

\ArticleName{Report from the Open Problems Session at OPSFA13\footnote{This paper is a~contribution to the Special Issue on Orthogonal Polynomials, Special Functions and Applications. The full collection is available at \href{http://www.emis.de/journals/SIGMA/OPSFA2015.html}{http://www.emis.de/journals/SIGMA/OPSFA2015.html}}}

\Author{Edited by Howard S.~COHL}

\AuthorNameForHeading{H.S.~Cohl}

\Address{Applied and Computational Mathematics Division, National Institute of Standards\\
 and Technology (NIST), Mission Viejo, CA, 92694 USA}
\Email{\href{mailto:howard.cohl@nist.gov}{howard.cohl@nist.gov}}
\URLaddress{\url{http://www.nist.gov/itl/math/msg/howard-s-cohl.cfm}}

\ArticleDates{Received January 28, 2016, in final form July 12, 2016; Published online July 21, 2016}

\Abstract{These are the open problems presented at the 13th International Symposium on Orthogonal Polynomials, Special Functions and Applications (OPSFA13), Gaithersburg, Maryland, on June 4, 2015.}


\Classification{31C12; 32Q10; 33C05; 33C45; 33C55; 33F10; 35J05}

\renewcommand{\thefootnote}{\arabic{footnote}}
\setcounter{footnote}{0}

\newcommand{\AuthorS}[1]{\noindent {\it #1} \par\vspace{1mm}\par}
\newcommand{\AddressS}[1]{\vspace{1mm}\par\noindent {\it #1} \par}

\section{Schur's inequality}
\AuthorS{Posed by Richard Askey}
\AddressS{Department of Mathematics, University of Wisconsin, Madison, WI, 53706, USA}
\Email{\href{mailto:askey@math.wisc.edu}{askey@math.wisc.edu}}

\bigskip

In Hardy, Littlewood and P\'olya's book ``Inequalities'' \cite[Problem 60 on p.~64]{HardyLittlewoodPolya}, the following inequality (communicated by I.~Schur) was stated:
\begin{gather*} x^n(x-y)(x-z)+y^n(y-x)(y-z)+z^n(z-x)(z-y) > 0, \end{gather*}
when $x$, $y$, $z$ are positive and not all equal, and $n \geq 0$. It is not hard to show that the inequality is true for all real $x$, $y$, $z$ for $n$ even when $>0$ is replaced by $\geq 0$.

There is a strange theorem of Hilbert~\cite{Hilbert1888}. He proved that if one has a polynomial in $k$ variables which is homogeneous of degree $j$ and which is nonnegative for all real values of the variables, then it can be written as a sum of $k$ squares when either $k$ or $j$ is $2$, and when $k=3$ and $j=4$, but not necessarily in all other cases. The case $n=2$ fits the exceptional condition $k=3$, $j=4$, so a natural question is: \textit{what is this representation?}

This problem has been solved in the mean time. The answer is
\begin{gather*} \frac14 \bigl( \big(2x^2-y^2-z^2+2yz-xz-xy\big)^2+3\big(y^2-z^2+xz-xy\big)^2 \bigr).
\end{gather*}
James Wan (Singapore University of Technology and Design) was the first to send this to me. Shortly after that, I got a solution using software of Erich Kaltofen (North Carolina State University) and Zhengfeng Yang (Shanghai key Laboratory of Trustworthy Computing), which was run by Zhengfeng Yang. Manuel Kauers (RISC, Johannes Kepler University of Linz, Austria) sent the problem to Erich Kaltofen, who was traveling and he sent it to his coworker. I was indeed hoping for something like this.

\section[Blumenthal--Nevai theorems for the one quarter class of orthogonal polynomials]{Blumenthal--Nevai theorems for the one quarter class\\ of orthogonal polynomials}

\AuthorS{Posed by Ted Chihara}
\AddressS{Purdue University, Calumet, Hammond, IN, 46323, USA}
\Email{\href{mailto:chihara@purduecal.edu}{chihara@purduecal.edu}}

\bigskip

Consider an orthogonal polynomial sequence (OPS), $\{P_n(x)\}$, which is defined by the classical three term recurrence relation ($n \ge 1$),
\begin{gather*}
P_n(x) = (x - c_n)P_{n-1}(x) - \lambda_n P_{n-2}(x),
\end{gather*}
with $P_{-1}(x)= 0$, $P_0(x) = 1$, $c_n$ real, and $\lambda_n > 0$. Suppose first that the sequences $\{c_n\}$ and~$\{\lambda_n\}$ converge to finite limits, $c$ and $\lambda$. Let $\sigma = \tfrac12(c - 2\sqrt{\lambda})$ and $\tau =\tfrac12(c + 2\sqrt{\lambda})$. In his dissertation written under David Hilbert, O.~Blumenthal~\cite{Blumenthal} proved that the zeroes of all of the~$P_n(x)$ form a set that is dense in the interval $[\sigma, \tau]$. More recently, P.~Nevai~\cite{Nevai1979} proved that in fact the interval $[\sigma, \tau]$ is the essential spectrum of the orthogonal polynomial system (i.e., of the operator given by the corresponding Jacobi matrix). See~\cite{Chihara68,Chihara03} for references and historical facts.

We would like to find analogues of these theorems when the sequences above are unbounded. To this end, let us consider the ``one quarter class'' of orthogonal polynomials (see \cite{Chihara91}); namely, let $\lim\limits_{n \to \infty}c_n = \infty$, and
\begin{gather*}
\lim_{n \to \infty} \frac{\lambda_{n+1}}{c_n c_{n+1}} = \frac{1}{4}.
\end{gather*}
Now let $x_{n,1} < x_{n,2} < \cdots < x_{n,n}$, denote the zeros of $P_{n}(x)$ and let
\begin{gather*}
\xi_i = \lim_{n \to \infty} x_{n,i}, \qquad \eta_j = \lim_{n \to \infty} x_{n,n-j+1}, \qquad
\sigma = \lim_{i \to \infty} \xi_i , \qquad \tau = \lim_{j \to \infty} \eta_j.
\end{gather*}
Now under the above conditions (i.e., for the one-quarter class), we will have $\tau = \infty$, and each of the three cases,
\begin{gather*}
\sigma = - \infty, \qquad |\sigma| < \infty, \qquad \sigma = \infty,
\end{gather*}
can occur. Sufficient conditions for each of the three cases to occur can be expressed in terms of the concept of ``eventual chain sequences'' (see \cite{Chihara82}).

Now the case $\sigma = \infty$ corresponds to the orthogonality measure having a discrete spectrum with $\infty$ as its only limit point. A specific example is furnished by certain Meixner polynomials of the first kind~\cite{chi1}. For the case when $|\sigma|<\infty$, we proved~\cite{Chihara68} that the set of all zeros of the corresponding orthogonal polynomials is dense in the interval $[\sigma, \infty)$. Later, we posed~\cite{Chihara03} the problem of determining, after imposing additional conditions, if necessary, that $[\sigma, \infty)$ is a~subset of the spectrum (that is, it is the essential spectrum). We hasten to remind that when~$\sigma$ is finite, the corresponding Hamburger moment problem will be determined (a~fact that was known to Stieltjes). Recently, Grzegorz \'{S}widerski~\cite{swiderski15a} has proven that, with a slight additional condition (a~certain monotonicity condition), in fact $[\sigma, \infty)$ is indeed the essential spectrum, thus settling this open problem (see also~\cite{AptekarevGeronimo14}). This now leaves the most difficult case, $\sigma=-\infty$.

We thus pose the inevitable question:~{\it Is there an analog of Blumenthal's theorem for this case?} The lack of specific examples of OPS whose orthogonality measure has a spectrum extending over the entire real line makes this situation the most difficult to conjecture about. The only known example of OPS of this type, are the Meixner polynomials of the second kind~\cite{chi1} (i.e., the so-called ``Meixner--Pollaczek'' polynomials). Do there exist any other OPS examples of this type?

\section[Generalized linearization formulas for hypergeometric orthogonal polynomials]{Generalized linearization formulas\\ for hypergeometric orthogonal polynomials}

\AuthorS{Posed by Howard S.~Cohl}
\AddressS{Applied and Computational Mathematics Division, National Institute of Standards\\
and Technology, Mission Viejo, CA, 92694, USA}
\Email{\href{mailto:howard.cohl@nist.gov}{howard.cohl@nist.gov}}

\bigskip

Given a hypergeometric orthogonal polynomial $P_n(x;{\bf a})$, where ${\bf a}$ is a set of arbitrary parameters, we would like to obtain closed-form expressions for the coefficients of generalized linearization formulas. A linearization formula for a hypergeometric orthogonal polynomial is an expression of type
\begin{gather}\label{linearization}
P_m(x;{\bf a})P_n(x;{\bf a})=\sum_{k=0}^{m+n} \alpha_{k,m,n}({\bf a}) P_k (x; {\bf a}) ,
\end{gather}
whereas a connection relation for a hypergeometric orthogonal polynomial is given by
\begin{gather}\label{connection}
P_n(x;{\bf a})=\sum_{k=0}^{n} \beta_{k,n}({\bf a},{\bf b}) P_k(x;{\bf b}).
\end{gather}
The coefficients $\alpha_{k,m,n}({\bf a})$, $\beta_{k,n}({\bf a},{\bf b})$ are usually given in terms of products of Pochhammer symbols (shifted factorials), generalized hypergeometric functions with fixed arguments, or multiple hypergeometric functions with fixed arguments. In our context, a generalized linearization formula is given by
\begin{gather*}
P_m(x;{\bf a})P_n(x;{\bf a}) =\sum_{k=0}^{m+n} \gamma_{k,m,n}({\bf a},{\bf b}) P_k (x; {\bf b}).
\end{gather*}
Generalized linearization formulae are obtained by inserting the connection relation (\ref{connection}) in the linearization formula (\ref{linearization}) and using series rearrangement with justification to identify the coefficient $\gamma_{k,m,n}({\bf a},{\bf b})$. Some concrete examples which we will discuss below are the Laguerre \cite[Section~9.12]{Koekoeketal}, Gegenbauer \cite[Section~9.8.1]{Koekoeketal}, continuous $q$-ultraspherical/Rogers \cite[Section~14.10.1]{Koekoeketal}, Jacobi \cite[Section~9.8]{Koekoeketal}, and continuous $q$-Jacobi \cite[Section~14.10]{Koekoeketal} polynomials.

Below we assume the product convention for ($q$-)Pochhammer symbols, namely $(a_1,a_2,\ldots,$ $a_n)_k:=(a_1)_k(a_2)_k\cdots(a_n)_k$, and $(a_1,a_2,\ldots,a_n;q)_k:=(a_1;q)_k(a_2;q)_k\cdots(a_n;q)_k$. We also adopt the following notations to indicate sequential elements in a list of elements, namely $\pm a:=\{a,-a\}$, and $f(g,\raisebox{-3pt}{$\choo{a}{b}$}):=\{f(g,a),f(g,b)\}$. For the linearization formulas that we present below, we will assume without loss of generality for $n,m\in {\mathbb N}_0:=\{0,1,2,\ldots\}$ that $n\ge m$.

For the Laguerre polynomials \cite[equation~(9.12.1)]{Koekoeketal}, the connection relation is \cite[equation~(18.18.18)]{NIST}
\begin{gather*}
L_n^{\alpha}(x)=\sum_{k=0}^n \frac{(\alpha-\beta)_{n-k}}{(n-k)!}L_k^\beta(x),
\end{gather*}
and the linearization formula is
\begin{gather*}
L_m^\alpha(x) L_n^\alpha(x)=\sum_{k=0}^{2m}A_{k,m,n}^{\alpha}L_k^\alpha(x),
\end{gather*}
where
\begin{gather*}
A_{k,m,n}^{\alpha}:=\frac{(-1)^m(-\alpha-n)_m(-2m,n-m+1)_k}{m!k!(\alpha+1+n-m)_k}\,\hyp32{-\alpha-m,\frac{-k+\choo01}{2}}{\,n-m+1,\frac12-m}{1}.
\end{gather*}
This linearization formula follows through a limiting procedure from the Rahman linearization formula of Jacobi polynomials given below,~(\ref{linJac}). Note that \cite[equations~(63) and (64)]{SanchezRuizetal} gives a~${}_3F_2(1)$ linearization formula for Laguerre polynomials, but that it does not seem to be correct. On the other hand, \cite[equations~(6.2) and (6.3)]{IsmailZeng2015} does give a correct linearization formula for Laguerre polynomials, but it is not easily converted to the relatively simple ${}_3F_2(1)$ form which results from Rahman's Jacobi linearization formula~(\ref{linJac}).

For the Gegenbauer polynomials \cite[equation~(9.8.19)]{Koekoeketal}, the connection relation is \cite[equation~(18.18.16)]{NIST}
\begin{gather*}
C_n^\lambda(x)= \sum_{k=0}^{\lfloor\frac{n}{2}\rfloor} \frac{\mu+n-2k}{\mu} \frac{(\lambda)_{n-k}(\lambda-\mu)_k}{k!(\mu+1)_{n-k}}C_{n-2k}^\mu(x),
\end{gather*}
and the linearization formula, assuming without loss of generality $n\ge m$, is (see \cite[equa\-tion~(18.18.22)]{NIST})
\begin{gather*}
C_m^\lambda(x) C_n^\lambda(x)= \sum_{k=0}^{m} B_{k,m,n}^{\lambda} C_{m+n-2k}^\lambda(x),
\end{gather*}
where
\begin{gather*}
B_{k,m,n}^{\lambda}:=\frac{(m+n+\lambda-2k)(m+n-2k)!(\lambda)_k(\lambda)_{m-k}(\lambda)_{n-k}(2\lambda)_{m+n-k}}{(m+n+\lambda-k)k!(m-k)!(n-k)!(\lambda)_{m+n-k}(2\lambda)_{m+n-2k}}.
\end{gather*}

For the continuous $q$-ultraspherical/Rogers polynomials \cite[equations~(14.10.17)]{Koekoeketal}, the connection relation is \cite[equations~(13.3.1)]{Ismail}
\begin{gather*}
C_n(x;\gamma|q)= \sum_{k=0}^{\lfloor\frac{n}{2}\rfloor} \frac{\beta^k(\gamma/\beta;q)_k(\gamma;q)_{n-k}(1-\beta q^{n-2k})} {(q;q)_k(\beta q;q)_{n-k}(1-\beta)} C_{n-2k}(x;\beta|q),
\end{gather*}
and the linearization formula is \cite[equation~(13.3.10)]{Ismail}
\begin{gather*}
C_m(x;\beta|q)C_n(x;\beta|q)= \sum_{k=0}^{m} D_{k,m,n}^{\beta} C_{m+n-2k}(x;\beta|q),
\end{gather*}
where
\begin{gather*}
D_{k,m,n}^{\beta}:=\frac{(q;q)_{m+n-2k}(\beta;q)_{m-k}(\beta;q)_{n-k}(\beta;q)_k(\beta^2;q)_{m+n-k}(1-\beta q^{m+n-2k})}{(q;q)_k(q;q)_{m-k}(q;q)_{n-k}(\beta q;q)_{m+n-k}(\beta^2;q)_{m+n-2k}(1-\beta)}.
\end{gather*}

The connection relation for Jacobi polynomials \cite[equation~(9.8.1)]{Koekoeketal} with two free parameters is given by (see for instance Ismail \cite[p.~256]{Ismail})
\begin{gather*}
P_n^{(\gamma,\delta)}(x)=\sum_{k=0}^n c_{k,n}^{\gamma,\delta;\alpha,\beta} P_k^{(\alpha,\beta)}(x),
\end{gather*}
where $\gamma,\delta>-1$, and such that if $\gamma,\delta\in(-1,0)$ then $\gamma+\delta+1\ne 0$,
\begin{gather*}
 c_{k,n}^{\gamma,\delta;\alpha,\beta}:=\frac{(\gamma+k+1)_{n-k}(n+\gamma+\delta+1)_k\Gamma(\alpha+\beta+k+1)}{(n-k)!\,\Gamma(\alpha+\beta+2k+1)} \\
\hphantom{c_{k,n}^{\gamma,\delta;\alpha,\beta}:=}{} \times\,{}_3F_2\left(\begin{array}{@{}c@{}}-n+k,n+k+\gamma+\delta+1,\alpha+k+1\\
\gamma+k+1,\alpha+\beta+2k+2\end{array};1\right).
\end{gather*}
It was originally thought that the linearization coefficients of the Jacobi polynomials were most simply represented by a double hypergeometric series \cite[p.~40]{AskeyRob2}, \cite[equation~(3.6) and (3.7)]{Miller68}. However, as pointed out to the author recently by Richard Askey, Rahman was able to prove that the linearization coefficients of Jacobi polynomials can be represented as a very well-poised~${}_9F_8(1)$~\cite{Rahman81a}. Note that there was a minor typo in Rahman's published result, where in the linearization coefficient, the term $(-\alpha-\beta-2m)$ should have been the Pochhammer symbol $(-\alpha-\beta-2m)_k$. Let $m,n\in{\mathbb N}_0$, and without loss of generality $n\ge m$. The (corrected) linearization formula for the Jacobi polynomials given in Rahman \cite[cf.~p.~919]{Rahman81a} is
\begin{gather}\label{linJac}
P_m^{(\alpha,\beta)}(x)P_n^{(\alpha,\beta)}(x)=\sum_{k=0}^{2m}h_{k,m,n}^{\alpha,\beta}P_{k+n-m}^{(\alpha,\beta)}(x),
\end{gather}
where
\begin{gather*}
h_{k,n,m}^{\alpha,\beta}:=\frac{(\alpha+1,\beta+1)_n(\alpha+\beta+1)_{2n-2m}(\alpha+\beta+1)_{2m}(\alpha+\beta+1+2n-2m+2k)}{m!(\alpha+\beta+1)_m(\alpha+1,\beta+1)_{n-m}(\alpha+\beta+2)_{2n}(\alpha+\beta+1)}\\
{}\times\frac{(n-m+1,\alpha+\beta+1+2n-2m,2\alpha+2\beta+2+2n,-2m,\alpha-\beta)_k}{k!(2\beta+2+2n-2m,\alpha+1+n-m,\alpha+\beta+2+2n,-\alpha-\beta-2m)_k} \\
{}\times\hyp98 { \beta\!+\!n\!-\!m\!+\!\frac12, \frac{\beta+n-m+\frac52}{2}, \beta\!+\!\frac12, \beta\!+\!n\!+\!1, -\alpha\!-\!m, \frac{\alpha+\beta+k+2n-2m+\choo{1}{2} }{2},
 \frac{-k+\choo{0}{1}}{2} } { \frac{\beta+n-m+\frac12}{2}, n\!-\!m\!+\!1, \frac12\!-\!m, \alpha\!+\!\beta\!+\!n\!+\!\frac{3}{2}, \frac{\beta-\alpha-k+\choo{1}{2}}{2}, \beta\!+\frac{k+2n-2m+\choo{2}{3}}{2}}{1}.
\end{gather*}

There is a corresponding result proved by Rahman \cite{Rahman81b} for the continuous $q$-Jacobi polynomials \cite[equation~(14.10.1)]{Koekoeketal} whose linearization coefficients are given in terms of a very well-poised ${}_{10}\phi_{9}(q,q)$, namely
\begin{gather*}
P_m^{(\alpha,\beta)}(x|q)P_n^{(\alpha,\beta)}(x|q)=\sum_{k=0}^{2m}g_{k,m,n}^{\alpha,\beta,q}P_{k+n-m}^{(\alpha,\beta)} (x|q),
\end{gather*}
where
\begin{gather*}
g_{k,m,n}^{\alpha,\beta,q}:=\frac{(q^{\alpha+1},q^{\beta+1};q)_n(q^{\alpha+\beta+\choo{2}{3}};q)_{n-m}(q^{\alpha+\beta+\choo{1}{2}};q)_{m}}
{(q^{\alpha+1},q^{\beta+1};q)_{n-m}(q^{\alpha+\beta+\choo{2}{3}};q)_{n}(q,q^{\alpha+\beta+1};q)_{m}}\\
\times\frac{(q^{n-m+1},q^{\frac{\alpha+\beta+3}{2}+n-m};q)_k(q^{-m},q^{\alpha+\beta+1+n},q^{\frac{\alpha+\beta+1}{2}+n-m},q^{\frac{\alpha-\beta}{2}};q^\frac12)_k}
{(q^{\alpha+1+n-m},q^{\frac{\alpha+\beta+1}{2}+n-m};q)_k(q^\frac12,q^{\frac{\alpha+\beta+2}{2}+n},q^{\beta+1+n-m},q^{-m-\frac{\alpha+\beta}{2}};q^\frac12)_k}
q^{\frac{m}{2}+\alpha(m-k)}\\
\times\qhyp{10}{9}{q^{\beta+n-m+\frac12},\pm q^{\frac{\beta+n-m+\raisebox{3pt}{\tiny $\frac52$}}{2}},q^{\beta+\frac12},q^{\beta+n+1},q^{-\alpha-m},
q^{n-m+\frac{\alpha+\beta+k+\choo12}{2}},q^{\frac{-k+\choo01}{2}}}{\pm q^{\frac{\beta+n-m+\raisebox{3pt}{\tiny $\frac12$}}{2}},q^{n-m+1},q^{\frac12-m},q^{\alpha+\beta+n+\frac32},
q^{\frac{\beta-\alpha-k+\choo12}{2}},q^{\beta+n-m+\frac{k+\choo23}{2}}}{q}{q}.
\end{gather*}

The best chance for finding generalized linearization coefficients which are hypergeometric functions is for the Gegenbauer and continuous $q$-ultraspherical/Rogers polynomials. This is because these linearization coefficients are given by products of Pochhammer symbols. Perhaps other orthogonal polynomials in the ($q$-)Askey scheme are amenable to this calculation, but we have yet to uncover further closed-form linearization formulae.

Motivations for considering ordinary linearization formulas and for connection formulas are given in~\cite{AskeyRob2}. Generalized linearization formulas, have the same motivations amplified by an ability to freely choose parameters. It has been suggested by an editor of the current special issue that the most simple example of a generalized linearization formula involves Chebyshev polynomials of the first kind \cite[equation~(6.4.13)]{AAR}
\begin{gather*}
T_n(x)= \frac{1}{\epsilon_n}\lim_{\mu\to 0}\frac{n+\mu}{\mu}C_n^\mu(x)=\frac{1}{\epsilon_n}\lim_{\mu\to 0}\frac{(\mu+1)_n}{(\mu)_n}C_n^\mu(x),
\end{gather*}
where $\epsilon_n=2-\delta_{n,0}$ is the Neumann factor, and Chebyshev polynomials of the second kind
\begin{gather*}
U_n(x)=C_n^1(x).
\end{gather*}
For generic values of $m$ and $n$, one has the following classical relations between the Chebyshev polynomials of the first and second kind. For instance, we have the following linearization formula \cite[equation~(5.1)]{AskeyRob2}
\begin{gather*}
T_m(x)T_n(x)=\frac12(T_{m+n}(x)+T_{m-n}(x)),
\end{gather*}
and interrelation formula \cite[equation~(18.9.9)]{NIST}
\begin{gather*}
T_n(x)=\frac12(U_n(x)-U_{n-2}(x)).
\end{gather*}
Hence there is
\begin{gather}
T_m(x)T_n(x)=\frac14\big(U_{m+n}(x)-U_{m+n-2}(x) +U_{m-n}(x)-U_{m-n-2}(x)\big).\label{productCheby}
\end{gather}

The first two formulas are rewritings of standard trigonometric identities. The product formula~(\ref{productCheby}) is a degenerate case of the generalized linearization formula for Gegenbauer polynomials
\begin{gather*}
C_m^\lambda(x)C_n^\lambda(x)=\sum_{k=0}^{\lfloor\frac{m+n}{2}\rfloor} f_{k,m,n}^{\lambda,\mu} C_{m+n-2k}^\mu(x),
\end{gather*}
where $f_{k,m,n}^{\lambda,\mu}\in{\mathbb R}$, for $\lambda,\mu\in(-1/2,\infty) \setminus\{0\}$. We are able to compute explicitly the coefficients $f_{k,m,n}^{\lambda,\mu}$ which contain a terminating, balanced, well-poised ${}_9F_8(1)$ which satisfies the transformation \cite[equation~(7.6.1)]{Bailey64}. We plan to publish this generalized linearization formula elsewhere.

\renewcommand{\thefootnote}{$*$}

\section{Simplifying multiple summations}

\AuthorS{Posed by Charles F.~Dunkl~$^a$ and Christoph Koutschan~$^b$\,\footnote{Christoph Koutschan was supported by the Austrian Science
Fund (FWF):~W1214.}}
\AddressS{$^{a)}$~Department of Mathematics, University of Virginia, Charlottesville VA, 22904, USA}
\EmailDD{\href{mailto:cfd5z@virginia.edu}{cfd5z@virginia.edu}}
\AddressS{$^{b)}$~Johann Radon Institute for Computational and Applied Mathematics (RICAM),\\
\hphantom{$^{b)}$}~Austrian Academy of Sciences, Linz, Austria}
\EmailDD{\href{mailto:christoph.koutschan@ricam.oeaw.ac.at}{christoph.koutschan@ricam.oeaw.ac.at}}

\bigskip

Here we describe a general problem area, rather than a specific open problem. In various computations, such as connection coefficients of families of orthogonal polynomials, or multiple iterated integrals, one arrives at a~multiple sum of hypergeometric form. The classical series of this type are the Lauricella series but typically more general series arise in practice (more parameters for example). What is needed is a systematic approach to find simplification to lower order summations when this is possible, or hopefully, when the answer is known. Roughly speaking, one would want a collection of known formulas, like the single-sum hypergeometric formulas with famous names (Gauss, Saalsch\"{u}tz, Dixon, Watson, \dots). The state-of-the-art today includes techniques for deriving recurrence or differential equations, by treating parameters in the sum as variables; these techniques are referred to as the holonomic systems approach~\cite{WilfZeilberger92,Zeilberger90}, and they are mainly based on the idea of creative telescoping. Several algorithms in this spirit have been proposed, for example Zeilberger's~\cite{Zeilberger90a}, Takayama's~\cite{Takayama90b}, and Chyzak's~\cite{Chyzak00}, to name just a few of them. These algorithms work especially well when there is a closed form (products of Pochhammer symbols, gamma functions, etc.), as the corresponding recurrence equation is of first order and can easily be solved. However, when the resulting recurrence is of higher order, it is more involved to find a nicer representation of the original sum, for example, as a single-sum. To some extent the algorithms developed by Schneider~\cite{Schneider05} in the framework of difference fields can be applied for this purpose. Nevertheless, it appears that often some human insight is needed. For example one may have to postulate the form of the single sum and then apply algorithms to prove the validity. We consider this as one of the key ingredients in this problem area. We illustrate these ideas with two worked-out examples.

The first example is a double sum, which comes in a terminating and in a non-terminating version; both have closed forms. For $m,n=0,1,2,\ldots$ we define
\begin{gather*}
S( m,n) :=\sum_{i=0}^{m}\frac{( -m) _{i}(n+1) _{i}}{i!( m+n+2) _{i}}\sum_{j=0}^{n}\frac{(-n) _{j}\left( \frac{1}{2}-n\right) _{j}}{j!\left( \frac{1}{2}\right) _{j}}\frac{1}{i+j+\frac{1}{2}}\\
\hphantom{S( m,n)}{} =2^{2m+2n}\frac{m!( m+n) !( m+n+1) !\left(\frac{1}{2}\right) _{n}}{n!( n+2m+1) !\left( \frac{1}{2}\right) _{m+n+1}}.
\end{gather*}
The sum is of double hypergeometric series form because
\begin{gather*}
 \frac{1}{i+j+\frac{1}{2}}=2\frac{\left(\frac{1}{2}\right)_{i+j}}{\left(\frac{3}{2}\right)_{i+j}}.
\end{gather*}
The sum is from~\cite{DunklGasper14}; the application is in~\cite{DunklSlater15}. There is a non-terminating form of this sum: for $n=0,1,2,\dots, \beta\notin\mathbb{N}$,
\begin{gather*}
S( \beta,n) =\sum_{j=0}^{\infty}\frac{( -\beta)
_{j}(n+1) _{j}}{j!( \beta+n+2) _{j}}\sum_{i=0}^{n}\binom{2n}{2i}\frac{1}{i+j+\frac{1}{2}}\\
\hphantom{S( \beta,n)}{} =2^{2\beta+2n}B\left( \beta+1,n+\frac{1}{2}\right) \frac{\Gamma(n+\beta+2) \Gamma( n+\beta+1) }{n!\Gamma(
n+2\beta+2) }.
\end{gather*}
This is so far unpublished; one proof relies on the Rogers--Dougall formula \cite[equation~(16.4.9)]{NIST}. Alternatively, one can evaluate this double sum by computer algebra algorithms. We first apply Zeilberger's algorithm~\cite{Zeilberger90a} to the inner sum (summation with respect to~$j$), as it is implemented in the \href{http://www.risc.jku.at/research/combinat/software/HolonomicFunctions}{\sf HolonomicFunctions} package~\cite{Koutschan09}, developed by one of us. It computes the second-order recurrence
\begin{gather*}
(i+n+2) (2i+2n+5) T(i+2,n) \\
\qquad{} = \big(4i^{2}+4in+12i+2n^{2}+5n+9\big) T(i+1,n) - (i+1) (2 i+1) T(i,n),
\end{gather*}
where
\begin{gather*}
T(i,n) := \sum_{j=0}^{n}\frac{(-n)_{j}\left( \frac{1}{2}-n\right) _{j}}{j!\left( \frac{1}{2}\right) _{j}}\frac{1}{i+j+\frac{1}{2}},
\end{gather*}
and a similar, again second-order, recurrence with respect to the parameter~$n$. We find that $T(i,n)$ is not a hypergeometric term, so
Zeilberger's algorithm cannot be applied to perform the summation with respect to~$i$. Instead, we use its generalization, Chyzak's algorithm~\cite{Chyzak00}. It computes the two first-order recurrences
\begin{gather*}
(n+1) (2m+n+2) (2m+2n+3) S(m,n+1)\\
\qquad{} = 4 (2n+1) (m+n+1) (m+n+2) S(m,n)
\end{gather*}
and
\begin{gather*}
(2m+n+2) (2m+n+3) (2m+2n+3) S(m+1,n)\\
\qquad {} = 8 (m+1) (m+n+1) (m+n+2) S(m,n),
\end{gather*}
from which the closed-form evaluation readily follows.

Our second example is the reduction of a double sum to a single sum (the problem arose in an integral over the compact group ${\rm Sp}(2)$~\cite{Dunkl07}). For $( \alpha_{1},\alpha_{2},\alpha_{3},\alpha_{4}) \in \mathbb{N}_{0}^{4}$ such that $\alpha_{1}\equiv\alpha_{2}\equiv\alpha _{3}\equiv\alpha_{4}\operatorname{mod}2$ (all even or all odd), let
\begin{gather*}
 b_{0}=\frac{\alpha_{2}+\alpha_{3}}{2},\qquad
 b_{1}=\frac{\alpha_{1}+\alpha_{4}}{2},\qquad
 b_{2}=\frac{\alpha_{2}+\alpha_{4}}{2},\qquad
 b_{3}=\frac{\alpha_{3}+\alpha_{4}}{2}.
\end{gather*}
Then consider the following double sum (a terminating Kamp\'{e} de F\'{e}riet double hypergeometric series of order~3 \cite[p.~244]{Erdelyi},
with~5 numerator and~3 denominator parameters, and argument~(1,1))
\begin{gather*}
s ( \alpha_{1},\alpha_{2},\alpha_{3},\alpha_{4} )
:=\frac{( 2\kappa) _{2b_{1}}( 2\kappa) _{2b_{0}
}\left( \frac{1}{2}\right) _{b_{1}}\left( \frac{1}{2}\right) _{b_{0}
}\left( \frac{1}{2}\right) _{b_{3}}}{( 4\kappa) _{2b_{1}
+2b_{0}}\left( \kappa+\frac{1}{2}\right) _{b_{1}}\left( \kappa+\frac{1}
{2}\right) _{b_{0}}\left( \kappa+\frac{1}{2}\right) _{b_{3}}}\\
\hphantom{s ( \alpha_{1},\alpha_{2},\alpha_{3},\alpha_{4} ):= }{}
 \times\sum_{i=0}^{\lfloor \alpha_{4}/2\rfloor }\sum
_{j=0}^{\lfloor \alpha_{3}/2\rfloor }\frac{( -\alpha_{4}) _{2i}( -\alpha_{3}) _{2j}( \kappa)
_{i+j}}{i!j!\left( \frac{1}{2}-b_{1}\right) _{i}\left( \frac{1}{2}
-b_{0}\right) _{j}\left( \frac{1}{2}-b_{3}\right) _{i+j}}2^{-2i-2j}\\
\hphantom{s ( \alpha_{1},\alpha_{2},\alpha_{3},\alpha_{4} )}{}
 =:\frac{( 2\kappa) _{2b_{1}}( 2\kappa) _{2b_{0}}
}{( 4\kappa) _{2b_{1}+2b_{0}}}s^{\prime}( \alpha_{1},\alpha_{2},\alpha_{3},\alpha_{4}) .
\end{gather*}
The simplification is
\begin{gather*}
s^{\prime}( \alpha_{1},\alpha_{2},\alpha_{3},\alpha_{4})
=\frac{\left( \frac{1}{2}\right) _{b_{1}}\left( \frac{1}{2}\right)
_{b_{2}}\left( \frac{1}{2}\right) _{b_{3}}}{\left( \kappa+\frac{1}{2}\right) _{b_{1}}\left( \kappa+\frac{1}{2}\right) _{b_{2}}\left(
\kappa+\frac{1}{2}\right) _{b_{3}}}\!
\sum_{i=0}^{\lfloor \alpha_{4}/2\rfloor }\!\!\frac{\left(
{-}\frac{\alpha_{4}}{2}\right) _{i}\left( \frac{1-\alpha_{4}}{2}\right)
_{i}( \kappa) _{i}( {-}\kappa-b_{1}-b_{0}) _{i}}{i!\left( \frac{1}{2}-b_{1}\right) _{i}\left( \frac{1}{2}-b_{2}\right)
_{i}\left( \frac{1}{2}-b_{3}\right) _{i}}.
\end{gather*}
The latter sum is a terminating balanced $_{4}F_{3}$-series. By use of the Whipple transformation it can be shown that $s^{\prime}$ is completely symmetric in its arguments. An intermediate step in formulating the single sum was to discover the recurrence
\begin{gather*}
\alpha_{1}\alpha_{4}\left(\kappa+\frac{1}{2}(\alpha_{2}+\alpha_{3}+1 )\right) s^{\prime}(\alpha_{1}-1,\alpha_{2}+1,\alpha_{3}+1,\alpha_{4}-1) \\
\qquad\quad{} +\frac{1}{2}\big( \alpha_{2}\alpha_{3}( \alpha_{1}+\alpha_{4}+1) -\alpha_{1}\alpha_{4}( \alpha_{2}+\alpha_{3}+1)
\big) s^{\prime}( \alpha_{1},\alpha_{2},\alpha_{3},\alpha_{4})\\
\qquad{} =\alpha_{2}\alpha_{3}\left( \kappa+\frac{1}{2}( \alpha_{1}+\alpha_{4}+1) \right) s^{\prime}( \alpha_{1}+1,\alpha_{2}-1,\alpha_{3}-1,\alpha_{4}+1) .
\end{gather*}

Once the single-sum representation is conjectured, it is again more or less routine to prove that it is equal to the original double sum $s(\alpha_{1},\alpha_{2},\alpha_{3},\alpha_{4})$, although the computations get a bit more involved now. Here it is convenient to consider the two cases (even and odd) separately. We found that Takayama's algorithm~\cite{Takayama90b} works best in this example, again using the implementation described in~\cite{Koutschan09}. For each side of the identity, it derives a set of three-term recurrence equations in the parameters~$\alpha_{1}$, $\alpha_{2}$, $\alpha_{3}$, $\alpha_{4}$; in more technical terms this means that the two annihilators have both holonomic rank~$2$. It turns out that the recurrences for the left-hand side perfectly agree with those for the right-hand side. Hence by comparing a few initial values the identity is established.

\section[Four open problems in orthogonal polynomials and random matrices]{Four open problems in orthogonal polynomials\\ and random matrices}

\AuthorS{Posed by Sheehan Olver\,\footnote{Sheehan Olver would like to acknowledge Deniz Bilman, Andrew Swan,
Alex Townsend, and Thomas Trogdon for helping to pose his list of questions.}}
\AddressS{School of Mathematics and Statistics, The University of Sydney, New South Wales, Australia}
\Email{\href{mailto:Sheehan.Olver@sydney.edu.au}{Sheehan.Olver@sydney.edu.au}}

\vspace{-1mm}

\subsection{Existence of a fast discrete spherical harmonic transform}

{\it Does there exist an algorithm that can efficiently convert from values of a function evaluated on a~grid on the sphere to spherical harmonic coefficients, and vice-versa?} Ideally, such an algorithm would be roughly of complexity $\mathcal{O}(n\log n)$ for a~grid of~$n$ points.

\vspace{-1mm}

\subsection{Easy-to-use software for uniform asymptotics of orthogonal polynomials}

\looseness=-1
{\it Is it possible to make uniform asymptotics with error bounds easy-to-use, for general orthogonal polynomials?} Recent work on quadrature~\cite{FastLegendreQuad,FastGaussQuad} and fast transforms~\cite{FastChebTransform} is built-up from uniform asymptotics, where error bounds are necessary to ensure accuracy and to optimize complexity. Riemann--Hilbert problems allow for uniform asymptotics for general orthogonal polynomials~\cite{DeiftStrongAsymptotics,DeiftOrthogonalPolynomials,OPsUI}, however, the methodology is hard to use for non-experts. A~software package that would take in a general weight and return the uniform asymptotic expansion would be ideal.

\vspace{-1mm}

\subsection{Spectrum of a finite-dimensional random symmetric Bernoulli matrix}

{\it What is the spectrum of a finite-dimensional random symmetric Bernoulli matrix?} A random symmetric Bernoulli matrix consists of entries that are randomly $\pm 1$, subject to a symmetry condition. What is the spectrum of a finite-dimensional random symmetric Bernoulli matrix? Unlike other symmetric random matrices such as Gaussian Orthogonal Ensemble (GOE), these only have a finite number of configurations of eigenvalues. While such matrices fall into the general framework of universality for Wigner ensembles \cite{YauUniversality,TaoUniversality} (describing the asymptotics of the spectrum), this does not explain the finite-dimensional picture.

\vspace{-1mm}

\subsection[Existence of a Wigner-like family corresponding to general invariant ensembles]{Existence of a Wigner-like family corresponding\\ to general invariant ensembles}

{\it Does there exist a Wigner-like family corresponding to general invariant ensembles?} For special cases we know such families exist. For example, the Wigner ensembles have the same limiting spectral density as a Gaussian Unitary Ensemble (GUE). Similarly, the Wishart ensembles have the same limiting spectral density as a Laguerre Unitary Ensembles (LUE).

\section{Positivity of an integral involving Gegenbauer polynomials}

\AuthorS{Posed by Rick Beatson~$^a$, Wolfgang zu Castell~$^b$ and Yuan Xu~$^c$}
\AddressS{$^{a)}$~Department of Mathematics and Statistics, University of Canterbury,\\
\hphantom{$^{a)}$}~Christchurch, New Zealand}
\EmailDD{\href{mailto:r.beatson@math.canterbury.ac.nz}{r.beatson@math.canterbury.ac.nz}}
\AddressS{$^{b)}$~Department of Scientific Computing, Helmholtz Zentrum M\"{u}nchen,\\
\hphantom{$^{b)}$}~German Research Center for Environmental Health, Neuherberg, Germany}
\EmailDD{\href{mailto:castell@helmholtz-muenchen.de}{castell@helmholtz-muenchen.de}}
\AddressS{$^{c)}$~Department of Mathematics, University of Oregon, Eugene, OR, 97403, USA}
\EmailDD{\href{mailto:yuan@uoregon.edu}{yuan@uoregon.edu}}

\bigskip

The following conjecture was stated in \cite[Conjecture 1.4]{BeatsonzuCastellXu2014}.

\begin{conj}\label{conj}
Let $\delta>0$, $\lambda > 0$ and $n\in{\mathbb N}_0$. For every $0<t<\pi$, define
\begin{gather*}%\label{eq:define_F}
F_n^{\lambda,\delta}(t) = \int_0^t (t-\theta)^\delta C_n^\lambda(\cos\theta) (\sin\theta)^{2\lambda} d\theta.
\end{gather*}
Then $F_n^{\lambda,\delta}(t)> 0$ for all $t$ in $(0,\pi]$ if and only if $\delta\geq\lambda+1$.
\end{conj}

It is known that if $F_n^{\lambda,\delta}(t) \ge 0$, then $F_n^{\lambda,\gamma}(t) \ge 0$ for $\gamma > \delta$. For
$\lambda > 0$, let $F_n^\lambda(t): = F_n^{\lambda, \lambda+1}(t)$. It is proved in \cite{BeatsonzuCastellXu2014} that $F_n^\lambda (t) \ge 0$
if $\lambda = \frac{d-2}{2}$ and $ d= 4,6,8$.

The conjecture is associated with the study of positive definite functions on the unit sphere. Under the assumption that $F_n^\lambda$ is nonnegative for $\lambda = \frac{d-2}{2}$ and all $n \in {\mathbb N}_0$, a P\'olya criterion for positive definite functions on the sphere ${\mathbb S}^{d-1}$ is established in~\cite{BeatsonzuCastellXu2014}.

\section[A family of polynomials related to a multiple zeta values identity]{A family of polynomials related\\ to a multiple zeta values identity}

\AuthorS{Posed by Wadim Zudilin}
\AddressS{School of Mathematical and Physical Sciences, University of Newcastle,\\ Callaghan NSW 2308, Australia}
\Email{\href{mailto:wzudilin@gmail.com}{wzudilin@gmail.com}}

\bigskip

The multiple zeta values (MZVs) are defined for positive integers $s_1,s_2,\dots,s_l$ with $s_1>1$ as the values of the convergent series
\begin{gather*}
\zeta(s_1,s_2,\dots,s_l) =\sum_{n_1>n_2>\dots>n_l\ge1}\frac1{n_1^{s_1}n_2^{s_2}\dotsb n_l^{s_l}}.
\end{gather*}
They satisfy numerous identities (some of which remain conjectural) and attract considerable interest of scientists working in number theory, algebraic geometry and mathematical physics.

The easiest identity is $\zeta(2,1)=\zeta(3)$. It was already given by Euler centuries ago. It in fact generalizes to the form
\begin{gather*}
\zeta(\underbrace{2,1,2,1,\dots,2,1}_{2l\text{ entries}})=\zeta(\underbrace{3,3,\dots,3}_{l\text{ entries}})
\end{gather*}
for $l=1,2,3,\dots$, an identity that can be proved by using a suitable integral representation of MZVs. There is a different proof of the identity in which certain biorthogonally looking polynomials show up. Namely, the polynomials form the one-parameter family
\begin{gather*}
B_n^\alpha(t) =\frac1{n!}\sum_{k=0}^n\frac{(\omega t)_k(\omega^2t)_k(\alpha+t)_{n-k}(\alpha-t+k)_{n-k}}{k!\,(n-k)!},
\end{gather*}
where $\omega=\exp(2\pi i/3)$ is the cubic root of unity. Though it is not obvious from the representation, we have $B_n^\alpha(t)\in\mathbb C[t^3]$ for $n=0,1,2,\dots$, so that we can view $B_n^\alpha$ as polynomials in $x=t^3$. To prove that $B_n^\alpha(t)\in{\mathbb C}[t^3]$, one can use the $3$-term recurrence relation
\begin{gather*}
%\label{rec-Ba}
\big((n+\alpha)^3-t^3\big)B_n^\alpha-(n+1)\big(2n^2+3n(\alpha+1)+\alpha^2+3\alpha+1\big)B_{n+1}^\alpha\\
\qquad{}+(n+2)^2(n+1)B_{n+2}^\alpha=0
\end{gather*}
and the initial conditions $B_0^\alpha=1$, $B_1^\alpha=\alpha^2$, satisfied by the polynomials. It is not hard to see from the recursion that the $x$-polynomials~$B_n^\alpha$ have degree~$[n/2]$. What is more surprising (but observed experimentally only) is that the zeroes of the polynomials are all real and follow a certain distribution on the negative half-line~$(-\infty,0)$.

A modified version of the MZV identity,
\begin{gather*}
\sum_{n_1>m_1>n_2>m_2>\dots>n_l>m_l\ge1} \frac{(-1)^{n_1+n_2+\dots+n_l}}{n_1^2m_1n_2^2m_2\dotsb n_l^2m_l}\\
\qquad{} =8^l\sum_{n_1>m_1>n_2>m_2>\dots>n_l>m_l\ge1} \frac{1}{n_1^2m_1n_2^2m_2\dotsb n_l^2m_l}
\end{gather*}
for $l=1,2,3,\dots$, has been recently established using a cumbersome machinery of MZVs. The identity is equivalent to proving that the polynomials $A_n(t)\in\mathbb Q[t^3]$ (of degree $[n/2]$ in $x=t^3$) produced by the recursion
\begin{gather*}
\big(n^3-(-1)^nt^3\big)A_n(t)+(n+1)^2(2n+1)A_{n+1}(t)+(n+2)^2(n+1)A_{n+2}(t)=0
\end{gather*}
and the initial conditions $A_0=1$, $A_1=0$ (no closed-form is known!) satisfy
\begin{gather*}
\sum_{k=0}^\infty A_k(t)=\prod_{j=1}^\infty\biggl(1+\frac{t^3}{8j^3}\biggr).
\end{gather*}
Equivalently, the polynomials $\widetilde A_n(t)=\sum\limits_{k=0}^nA_k(t)$ that come from the recursion
\begin{gather}\label{endrec}
\big(n^3-(-1)^nt^3\big)\widetilde A_{n-1}+(2n+1)n\widetilde A_n-(n+1)^2n\widetilde A_{n+1}=0
\end{gather}
satisfy
\begin{gather*}
\lim_{n\to\infty}\widetilde A_n(t)=\prod_{j=1}^\infty\left(1+\frac{t^3}{8j^3}\right).
\end{gather*}
Note that the zeroes of the polynomials $A_n(t)$ and $\widetilde A_n(t)$ are also expected to lie on the negative half-line $(-\infty,0)$. The details of the story can be found in~\cite{Zudilin15}.

\begin{Questions}
 Is there an argument to deduce the limit of $\widetilde A_n(t)$ as $n\to\infty$ using the recurrence relation for the polynomials~\eqref{endrec}? Can the polynomials $A_n(t)$ and $\widetilde A_n(t)$ be given in an explicit hypergeometric form?
\end{Questions}

\pdfbookmark[1]{References}{ref}
\LastPageEnding

\end{document}